\documentclass[oneside,notitlepage,12pt]{article}

\usepackage{amssymb}
\usepackage[leqno]{amsmath}
\usepackage{amsfonts}
\usepackage{amsopn}
\usepackage{amstext}
\usepackage{amsthm}

\usepackage{enumitem}

\usepackage{tikz}
\usetikzlibrary{cd}

\usepackage{verbatim}
\usepackage[colorlinks, backref]{hyperref}

\usepackage{calrsfs}
\usepackage{fourier-orns}
\usepackage{clock} 
\usepackage[alpine, weather]{ifsym}

\renewcommand{\Bbb}{\mathbb}

\newenvironment{pf}{\begin{proof}}{\end{proof}}

\newcommand{\Aa}{{\Bbb{A}}}

\newcommand{\Be}{{\Bbb{B}}}

\newcommand{\Qyu}{{\Bbb{Q}}}
\newcommand{\Err}{{\Bbb{R}}}

\newcommand{\al}{\alpha}

\newcommand{\sig}{\sigma}
\newcommand{\eps}{\varepsilon}
\renewcommand{\phi}{\varphi}
\renewcommand{\rho}{\varrho}

\newcommand{\rest}{\restriction}

\newcommand{\ntr}{{n\in\omega}}

\newcommand{\loe}{\leq}
\newcommand{\goe}{\geq}

\newcommand{\subs}{\subseteq}
\newcommand{\sups}{\supseteq}
\newcommand{\nnempty}{\ne\emptyset}

\newcommand{\ovr}{\overline}

\newcommand{\id}[1]{{\operatorname{i\!d}_{#1}}} 

\newcommand{\liminv}{\varprojlim}

\newcommand{\oraz}{\qquad\text{and}\qquad}

\newcommand{\meet}{\wedge}

\newcommand{\join}{\vee}

\newtheorem{result}{Theorem}
\newtheorem{tw}{Theorem}[section]
\newtheorem{lm}[tw]{Lemma}
\newtheorem{prop}[tw]{Proposition}

\theoremstyle{definition}

\theoremstyle{remark}

\newcommand{\setof}[2]{\{#1\colon #2\}}

\newcommand{\sett}[2]{\{#1\}_{#2}}
\newcommand{\sn}[1]{\{#1\}} 
\newcommand{\dn}[2]{\{#1,#2\}} 
\newcommand{\pair}[2]{\langle #1, #2 \rangle} 
\newcommand{\triple}[3]{\langle #1, #2, #3 \rangle} 
\newcommand{\map}[3]{#1\colon #2 \to #3} 
\newcommand{\img}[2]{#1[#2]} 
\newcommand{\inv}[2]{{#1}^{-1}[#2]} 

\newcommand{\cantor}{2^\omega}

\newcommand{\fra}{Fra\"iss\'e}

\newcommand{\bfK}{{\mathbf K}}

\providecommand{\nat}{\omega}

\newcommand{\obj}[1]{\operatorname{Obj}\left(#1\right)}

\newcommand{\fK}{{\mathfrak{K}}}
\newcommand{\fL}{{\mathfrak{L}}}

\newcommand{\fC}{{\mathfrak{C}}}

\newcommand{\cmp}{\circ} 

\newcommand{\separator}{\begin{center} \leafright \leafright \leafright \decotwo \decotwo \decotwo \leafleft \leafleft\leafleft
\end{center}}

\providecommand{\ar}{\arrow}

\newcommand{\Prob}{\operatorname{\mathbf P}}
\newcommand{\Probzero}{\operatorname{\mathbf P_0}}

\newcommand{\prob}{\operatorname{P}}
\newcommand{\los}{\ensuremath{\boxminus}}
\newcommand{\losv}[1]{\ensuremath{\boxminus_{#1}}}

\newcommand{\kfP}{{\mathfrak{P}_{\rm fin}}}
\newcommand{\kP}{{\mathfrak{P}_{\rm prof}}}
\newcommand{\kfPv}[1]{{\mathfrak{P}^{#1}_{\rm fin}}}
\newcommand{\kPv}[1]{{\mathfrak{P}^{#1}_{\rm prof}}}

\newcommand{\Jeden}{{\mathbf 1}}

\newcommand{\koment}[1]{}

\title{Homogeneous probability measures on the Cantor set}
\author{
Wojciech Bielas \\
{\small Institute of Mathematics, Czech Academy of Sciences, Czechia}\\
{\small Institute of Mathematics, University of Silesia in Katowice, Poland}
\and
Wies{\l}aw Kubi\'s \\
{\small Institute of Mathematics, Czech Academy of Sciences, Czechia}\\
{\small Institute of Mathematics, Cardinal Stefan Wyszy\'nski University in Warsaw, Poland}
\and
Marta Walczy\'nska \\
{\small Institute of Mathematics, Czech Academy of Sciences, Czechia}\\
{\small Institute of Mathematics, University of Silesia in Katowice, Poland}
}

\date{\clocktime\today}

\begin{document}

\maketitle

\begin{abstract}
We show that every homeomorphism between closed measure zero subsets extends to a measure preserving auto-homeomorphism, whenever the Cantor set is endowed with a suitable probability measure.
This is valid both for the standard product measure, as well as for the universal homogeneous rational measure.


\end{abstract}

\tableofcontents

\section{Introduction}

Our objects of study are profinite spaces with probability measures.
A topological space is \emph{profinite} if it is homeomorphic to the inverse limit of a sequence of nonempty finite sets.
It is well known that a nonempty topological space is profinite if and only if it is compact Hausdorff and has a countable basis consisting of clopen (i.e. closed and open at the same time) sets.
Moreover, it is known that every finitely additive measure on the algebra of all clopen sets extends uniquely to a Borel measure (see~\ref{SectExtRozMiry} for more information).

A measure $\mu$ on a profinite space $K$ is \emph{rational} (resp. \emph{strictly positive}) if $\mu U$ is a rational number (resp. $\mu U >0$) for every nonempty clopen set $U \subs K$.
Given two measure spaces $(X, \mu)$, $(Y, \nu)$, a mapping $\map f X Y$ is \emph{measure preserving} if $\mu \inv f B = \nu B$ for every measurable set $B \subs Y$. We shall use this concept for continuous surjections $f$ and regular Borel measures $\mu$, $\nu$ only.
Recall that every finite Borel measure on a compact metric space is regular.

Kechris \& Rosendal~\cite{KecRos} showed that the countable atomless Boolean algebra equipped with a homogeneous finitely additive probability measure taking positive rational values is a \fra\ limit. This result can be formalized in topological terms as follows:

\begin{result}\label{THMmejnOne}
	There exists a unique, up to measure preserving homeomorphisms, strictly positive rational probability measure $\Prob$ on the Cantor set $\cantor$, satisfying the following condition.
	\begin{enumerate}
		\item[{\rm(\los)}] For every nonempty clopen set $U \subs \cantor$, for every rational numbers $r_0, r_1 > 0$ with $r_0 + r_1 = \Prob U$, there exists a partition $\dn{U_0}{U_1}$ of $U$ into clopen sets such that $\Prob U_i = r_i$ for $i=0,1$.
	\end{enumerate}
	Furthermore, $\Prob$ is isomorphic to the product of all finite strictly positive rational probability measures.
\end{result}

We present a self-contained proof of this result in Section~\ref{SecPfThmOneJdn} below. In fact, Kechris and Rosendal~\cite{KecRos} showed that the Boolean algebra generated by all rational intervals contained in $[0,1]$ is a model for this structure, which can be clearly seen by applying condition (\los).

In the next results we assume that the Cantor set is endowed with the probability measure $\Prob$ from Theorem~\ref{THMmejnOne}.

\begin{result}\label{THMmejnTwo}
	For every $n>0$, for every positive rational numbers $r_0, \dots, r_{n-1}$ with $r_0 + \dots + r_{n-1} = 1$ there exists a partition $\{U_0, \dots, U_{n-1}\}$ into clopen subsets of $\cantor$ such that $\Prob U_i = r_i$ for $i<n$.
	Furthermore, if $\{V_0, \dots, V_{n-1}\}$ is another clopen partition with this property and $K \subs \cantor$ is a closed measure zero set such that $U_i \cap K = V_i \cap K$ for $i < n$, then there exists a measure preserving homeomorphism $\map h \cantor \cantor$ satisfying $\img h {U_i} = V_i$ for $i<n$ and $h \rest K = \id K$.
\end{result}

\begin{result}\label{THMmejnThree}
	Let $K$ be a profinite space with a strictly positive rational probability measure.
	Then there exist a topological embedding $\map \eta K \cantor$ and a measure preserving continuous surjection $\map f \cantor K$ such that $\Prob(\img \eta K) = 0$ and $f \cmp \eta = \id K$.
\end{result}

In fact, we shall prove more. Namely, if $K \subs \cantor$ is a closed measure zero set, then for \emph{every} rational strictly positive probability measure on $K$ there exists a measure preserving retraction $\map r \cantor K$.

\begin{result}\label{THMmejnFour}
	Let $\map h K L$ be a homeomorphism between closed measure zero subsets of $\cantor$. Then there exists a measure preserving homeomorphism $\map{H}{\cantor}{\cantor}$ satisfying $H \rest K = h$.
\end{result}

The last result is a measure-theoretic variant of the theorem of Knaster \& Reichbach~\cite{KnaRei} stating that homeomorphisms between closed nowhere dense subsets of $\cantor$ extend to auto-homeomorphisms of $\cantor$.
Theorem~\ref{THMmejnThree} is a strengthening of the well known fact saying that every profinite space is homeomorphic to a retract of the Cantor set. Note that every profinite space $K$ has a strictly positive rational probability measure, see Proposition~\ref{PROPsdbAIFHU} below.

In the last section we discuss possible variants of the results above, where the value set of the measure is a fixed countable subset of $[0,1]$, for example, consisting of dyadic rationals.

\section{Preliminaries}

We shall need some basic concepts from category theory. We will also use some results from~\cite{KubFra}. For undefined notions concerning category theory we refer to~\cite{MacLane}.

A category $\fC$ will be identified with its class of arrows; $\obj{\fC}$ will denote its class of objects.
Given $X,Y \in \obj{\fC}$ the set of $\fC$-arrows from $X$ to $Y$ will be denoted by $\fC(X,Y)$.
If $f \in \fC$ then we write $\map f X Y$ instead of $f \in \fC(X,Y)$, that is, $X$ is the domain of $f$ and $Y$ is the codomain of $f$.
An arrow $f \in \fC$ is called an \emph{epi} (often called an \emph{epimorphism}) if for every $\fC$-arrows $g_0,g_1$ the equation $g_0 \cmp f = g_1 \cmp f$ implies $g_0 = g_1$ (when writing expressions like $f \cmp g$, we implicitly assume that the domain of $f$ coincides with the codomain of $g$).

An \emph{inverse sequence} in $\fC$ is a contravariant functor from the category of natural numbers $\pair \omega \loe$ into $\fC$.
Specifically, an inverse sequence consists of a sequence $X_0, X_1, \ldots$ of $\fC$-objects together with a collection of $\fC$-arrows $\sett{x_n^m}{n \loe m < \nat}$
such that $\map{x_n^m}{X_m}{X_n}$ for $n \loe m$ and the following compatibility conditions are satisfied.
$$x_n^n = \id{X_n} \oraz (\forall\;k<\ell<m)\;\; x_k^m = x_k^\ell \cmp x_\ell^m.$$
Such a sequence will be denoted shortly by $\vec X$ and $x_n^m$ are called the \emph{bonding arrows} of $\vec X$.
The \emph{limit} of $\vec X$ is formally a pair $\pair{L}{\sett{x_n^\infty}{\ntr}}$, where $L \in \obj{\fC}$, $\map{x_n^\infty}{L}{X_n}$ satisfy $x_n^m \cmp x_m^\infty = x_n^\infty$ for $n<m$, and for every other pair $\pair{Y}{\sett{p_n}{\ntr}}$ satisfying $\map{p_n}{Y}{X_n}$, $x_n^m \cmp p_m = p_n$ for $n<m$, there exists a unique $\fC$-arrow $\map g Y L$ such that $x_n^\infty \cmp g = p_n$ for $\ntr$.
The limit (if it exists) is determined uniquely up to isomorphisms. The arrows $x_n^\infty$ are called the \emph{canonical projections} from the limit.
In the category of sets or topological spaces the limit of a sequence $\vec X$ can be chosen to be a suitable subset of the product $\prod_{\ntr}X_n$, and the canonical projections may be chosen to be restrictions of the usual projections onto the coordinates.

\subsection{Inverse \fra\ theory}

We recall the setup and basic results concerning \fra\ theory in the context of inverse sequences. The case of concrete categories of finite and profinite models was treated in~\cite{IrwSol}, while the abstract category-theoretic approach goes back to Droste \& G\"obel~\cite{DroGoe, DroGoe2}.
In fact, our setup fits into the framework of Droste \& G\"obel theory~\cite{DroGoe}.

Let $\fK$ be a fixed category.
We say that $\fK$ has the \emph{inverse amalgamation property} if for every $\fK$-arrows $f,g$ with the same codomain there exist $\fK$-arrows such that $f \cmp f' = g \cmp g'$ (in particular, $f',g'$ have the same domains).
We say that $\fK$ is \emph{inverse directed} if for every $X,Y \in \obj{\fK}$ there are $Z \in \obj{\fK}$ and $\fK$-arrows $\map f Z X$, $\map g Z Y$.
Finally, we say that $\fK$ \emph{is essentially countable} if $\fK$ has countably many objects up to isomorphisms, and for every $X,Y \in \obj{\fK}$ the hom-set $\fK(X,Y)$ is countable.
We shall say that $\fK$ is an \emph{inverse \fra\ category} if it is essentially countable, inverse directed, and has the inverse amalgamation property.

Let us now assume that $\fK$ is a full subcategory of a bigger category $\fL$ such that the following compatibility conditions are satisfied.
\begin{enumerate}[itemsep=0pt]
	\item[(L0)] All $\fL$-arrows are epi.
	\item[(L1)] Every inverse sequence in $\fK$ has the limit in $\fL$.
	\item[(L2)] Every $\fL$-object is the limit of an inverse sequence in $\fK$.
	\item[(L3)] For every inverse sequence $\vec X$ in $\fK$ with $K = \lim \vec X$ in $\fL$, for every $\fK$-object $Y$, for every $\fL$-arrow $\map f K Y$ there exist $n$ and a $\fK$-arrow $\map {f'}{X_n}Y$ such that $f = f' \cmp x_n^\infty$.
\end{enumerate}
Now, an $\fL$-object $U$ will be called \emph{$\fK$-generic} if
\begin{enumerate}[itemsep=0pt]
	\item[(G1)] $\fL(U,X) \nnempty$ for every $X \in \obj{\fK}$.
	\item[(G2)] For every $\fK$-arrow $\map f Y X$, for every $\fL$-arrow $\map g U X$ there exists an $\fL$-arrow $\map h U Y$ such that $f \cmp h = g$.
\end{enumerate}
An object satisfying (G2) will be called \emph{$\fK$-projective}. An object satisfying (G1) is sometimes called \emph{projectively $\fK$-universal}.
A standard argument shows that if $\fK$ has a weakly terminal object then (G1) is a consequence of (G2).
Indeed, if $X \in \obj{\fK}$ and $\map f X T$ is a $\fK$-arrow such that $T$ is weakly terminal, then by (L2) there is an $\fL$-arrow $\map j U T$.
Using (G2) we find an $\fL$-arrow $\map h U X$ satisfying $j = f \cmp h$, showing that $\fL(U,X) \nnempty$.

The following facts are consequences of~\cite[Cor. 3.8]{KubFra} applied to the opposite category of $\fK$.

\begin{tw}\label{ThmFraisseReversedOne}
	Assume $\fK \subs \fL$ satisfy {\rm(L0)}--{\rm(L3)} and $\fK$ is an inverse \fra\ category.
	Then there exists a unique, up to isomorphisms, $\fK$-generic object in $\fL$.
\end{tw}

\begin{tw}\label{ThmFraisseReversedTwo}
	Assume $\fK \subs \fL$ satisfy {\rm(L0)}--{\rm(L3)} and $U \in \obj{\fL}$ is $\fK$-generic.
	Then
	\begin{enumerate}[itemsep=0pt]
		\item[{\rm(1)}] $\fL(U, X) \nnempty$ for every $X \in \obj{\fL}$.
		\item[{\rm(2)}] For every $A \in \obj{\fK}$, for every $\fL$-arrows $\map {f_i}U A$, $i=0,1$, there exists an automorphism $\map h U U$ such that $f_1 = f_0 \cmp h$.
	\end{enumerate}
\end{tw}

Property (1) is called \emph{projective universality} while property (2) is called \emph{projective homogeneity} (see~\cite{IrwSol}).
Finally, let us note the converse to Theorem~\ref{ThmFraisseReversedOne}.

\begin{tw}\label{ThmFraisseRevConvers}
	Assume $\fK \subs \fL$ satisfy {\rm(L0)}--{\rm(L3)}.
	If a $\fK$-generic object exists in $\fL$, then $\fK$ is an inverse \fra\ category.	
\end{tw}

\begin{pf}
	Only the inverse amalgamation property requires an argument, however this can be easily deduced from (G1), (G2) and (L2), (L3), (L0).
	Specifically, if $U$ is $\fK$-generic and $\map f X Z$, $\map g Y Z$ are $\fK$-arrows then by (G1) we find an $\fL$-arrow $\map p U X$ and using (G2) we find an $\fL$-arrow $\map q U Y$ satisfying $f \cmp p = g \cmp q$.
	Using (L2) followed by (L3) we conclude that $p = p' \cmp u_n^\infty$ and $q = q' \cmp u_n^\infty$ for some $u_n^\infty$ coming from a fixed inverse sequence with limit $U$.
	Finally, by (L0) we obtain $f \cmp p' = g \cmp q'$.
\end{pf}

\subsection{Extending finitely additive measures}\label{SectExtRozMiry}

The problem of extending finitely additive measures to $\sig$-additive measures belongs to the folklore. For example, it had been extensively studied by Bachman \& Sultan~\cite{BaSu}, from which our special case can be extracted.  
Anyway, for the sake of completeness, we provide a direct argument below.
Actually, the metrizable case can be proved more directly by using the classical Carath\'eodory's extension theorem\footnote{The Referee pointed out to us that the Hahn-Kolmogorov theorem (cf. \url{https://en.wikipedia.org/wiki/Hahn-Kolmogorov_theorem}) could be applied, however it is a special case of the more well known Carath\'eodory's theorem.},
because the algebra of Borel sets is generated by the subalgebra of clopen sets and every measure $\mu$ defined on clopen subsets of a compact space is conditionally $\sig$-additive, i.e., satisfies $\mu(\bigcup_{\ntr}A_n) = \sum_{\ntr}\mu(A_n)$ whenever $A_n$ are pairwise disjoint clopen sets with $\bigcup_{\ntr} A_n$ clopen, since by compactness, $A_n = \emptyset$ for all but finitely many $\ntr$.

\begin{prop}\label{PropExtRozMiry}
	Let $K$ be a 0-dimensional compact space and let $\mu$ be a finitely additive measure defined on the algebra of all clopen subsets of $K$ whose values are in $[0, +\infty)$.
	Then there is a unique regular Borel $\sig$-additive measure $\ovr \mu$ extending $\mu$.
\end{prop}

\begin{pf}
	Given a continuous real-valued function $f \in C(K)$, define $\phi(f)$ to be the Riemann integral of $f$ with respect to $\mu$.
	More precisely, $\phi(f)$ is the supremum of all expressions of the form
	$$\sum_{i=1}^{n}f(x_i) \mu(U_i),$$
	where $\{U_1, \dots, U_n\}$ is a partition of $K$ into clopen sets, and $x_i \in U_i$ for each $i \loe n$.
	Following the classical theory of Riemann integrals, it is easy to see that $\phi$ is a positive bounded linear functional on $C(K)$. By the Riesz representation theorem, there exists a unique regular Borel measure $\ovr \mu$ such that $\phi(f) = \int_K f \,d\mu$ for every $f \in C(K)$.
	It is clear that $\ovr \mu$ extends $\mu$.
	Uniqueness follows from the fact that every regular Borel measure extending $\mu$ induces a linear functional as above and this functional must be equal to $\phi$ (coincidence of the Lebesgue integral with the Riemann integral on the space of continuous functions).
\end{pf}

It is well known that if $K$ is a metrizable compact space then every finite Borel measure is automatically regular.
Recall that a measure $\mu$ is \emph{regular at} a Borel set $S \subs K$ if for every $\eps>0$ there are a compact set $F \subs S$ and an open set $U \sups S$ with $\mu(U \setminus F) < \eps$. A measure is \emph{regular} if it is regular at every Borel set.
The family of all Borel sets $S$ such that $\mu$ is regular at $S$ is easily seen to be a $\sig$-algebra. Finally, if $K$ is compact metric then every open set is a countable union of compact subsets, therefore this $\sig$-algebra contains all open sets and hence coincides with the Borel $\sig$-algebra.

\section{Probability spaces as a category}

A finite probability space $\pair S {\prob_S}$ serves as a model for some real-life events, actions, behaviors, etc. The elements of $S$ are called \emph{elementary events}.
In practice, elementary events could be called indivisible events or events that from observer's point of view cannot be decomposed into sub-events.
In this context, a measure preserving surjection $\map{f}{\pair T{\prob_T}}{\pair S{\prob_S}}$ may represent an improvement of the probability space, in the sense that some of the elementary events are no longer elementary. In practice, this could be, for example, upgrading the equipment or the technology while observing certain events or performing a circular experiment.
Obviously, if $f$ is as above, then the space $\pair T{\prob_T}$ carries more information than $\pair{S}{\prob_S}$.
One can also argue that it suffices to restrict attention to rational probabilities, as in practice it is typically difficult or even impossible to compute the precise probabilities of events, therefore rational approximations indeed make sense.

We shall work with the category $\kfP$ whose objects are finite 
rational probability spaces $\pair S \prob$ such that $\prob \sn x > 0$ for every $x \in S$, that is, $\prob$ is strictly positive.  The arrows are measure preserving surjections.
Note that every probability space is nonempty, because it is required that the measure of the entire space is $1$, while the measure of the empty set is always $0$.

We shall also consider the category $\kP$, whose objects are profinite rational probability spaces $\pair K \prob$ whose probability is strictly positive. That is, $K \nnempty$ is compact Hausdorff, has a countable basis consisting of clopen sets, and $\prob U$ is a positive rational number for every nonempty clopen set $U \subs K$.
The arrows of $\kP$ are continuous measure preserving surjections.
Clearly, $\kfP$ is a full subcategory of $\kP$.
Note that $\kfP$ has a terminal object 
$\Jeden = \pair{\sn0}{P}$ that is also terminal in $\kP$. That is, for every $K \in \obj{\kP}$ there is a unique arrow $\map{1_K}{K}{\Jeden}$ (recall that all the $\kP$-objects are assumed to be nonempty).

A dual approach, using Boolean algebras and measure preserving embeddings, can be found in Kechris \& Rosendal~\cite{KecRos}.

\begin{lm}\label{LMdwaAjedn}
	Let
	$$\begin{tikzcd}
	\pair{S_0}{\prob_0} & \pair{S_1}{\prob_1} \ar[l, two heads, "s_0^1"'] & \cdots \ar[l, two heads, "s_1^2"']
	\end{tikzcd}$$
	be an inverse sequence in $\kfP$ and let $K$ be its inverse limit in the category of topological spaces.
	Then there exists a unique probability measure $\prob$ on $K$ such that all canonical projections $\map{s_n^\infty}{K}{S_n}$ become measure preserving.
	Furthermore, $\prob$ is rational and strictly positive.
\end{lm}

\begin{pf}
	Fix a clopen set $U \subs K$.
	Then $U = \inv{(s_n^\infty)}{V}$ for a unique $V \subs S_n$.
	Define $\prob U := \prob_n V$.
	Note that this does not depend on the choice of $n$, because the bonding mappings are measure preserving.
	If $U_0, U_1$ are pairwise disjoint clopen sets then we can find $m$ such that $U_i = \inv{(s_m^\infty)}{V_i}$ and $V_0 \cap V_1 = \emptyset$.
	Thus $\prob (U_0 \cup U_1) = \prob_m(V_0 \cup V_1) = \prob_m V_0 + \prob_m V_1 = \prob U_0 + \prob U_1$, showing that $\prob$ is finitely additive.
	Clearly, $\prob$ is strictly positive, because each $\prob_n$ is so.
	By Proposition~\ref{PropExtRozMiry}, $\prob$ extends uniquely to a Borel measure on $K$.
	Uniqueness of $\prob$ is obvious, as there is only one choice to define it on clopen sets.
\end{pf}

The following fact is trivial.

\begin{lm}\label{LMciongimiarp}
	Let $K$ be the inverse limit of a sequence $\vec S$ of finite nonempty sets and assume $\prob$ is a strictly positive rational probability measure on $K$.
	Then each $S_n$ has a rational strictly positive probability measure $\prob_n$, defined by $\prob_n A := \prob \inv{(s_n^\infty)}{A}$, so that all the bonding mappings (as well as the canonical projections $s_n^\infty$)
	become measure preserving.
\end{lm}

\begin{prop}\label{PROPsdbAIFHU}
	Every profinite space admits a strictly positive rational probability measure.
\end{prop}

\begin{pf}
	Let $K = \liminv \vec X$, where $\vec X$ is an inverse sequence of finite nonempty sets $X_n$, where the bonding mappings are surjections. It is clear how to define inductively rational strictly positive measures on each $X_n$ so that the bonding mappings become measure preserving. By Lemma~\ref{LMdwaAjedn}, this induces a strictly positive rational probability measure on $K$.
\end{pf}

\begin{lm}\label{LMomegaSmlns}
	Let $\vec S$ be an inverse sequence in $\kfP$ with limit $K \in \obj{\kP}$.
	Let $\map q K T$ be a $\kP$-arrow with $T$ finite.
	Then there exist $n$ and a $\kfP$-arrow $\map {p}{S_n} T$ such that $q = p \cmp s_n^\infty$.
\end{lm}

\begin{pf}
	As $q$ is continuous, there is $n$ such that $q = p \cmp s_n^\infty$ for some mapping $\map{p}{S_n}{T}$.
	Given $A \subs T$ we have
	$$\prob A = \prob \inv q A = \prob(\inv{(s_n^\infty)}{\inv{p}{A}}) = \prob_n(\inv p A),$$
	therefore $p$ is measure preserving, i.e., $p$ is a $\kfP$-arrow.
\end{pf}

\subsection{Pullbacks}

Let $\map f X Z$, $\map g Y Z$ be two mappings between nonempty sets.
Define
$$W = \setof{\pair x y \in X \times Y}{f(x) = g(y)},$$
and let $\map{\pi_f}{W}{X}$, $\map{\pi_g}{W}{Y}$ be the canonical projections, i.e., $\pi_f(x,y) = x$, $\pi_g(x,y) = y$ for $\pair xy \in W$.
The triple $\triple{W}{\pi_f}{\pi_g}$ will be called the \emph{pullback} of $\pair f g$.
In category theory, the commutative square
$$\begin{tikzcd}
Y \ar[d, "g"'] & W \ar[l, "\pi_g"'] \ar[d, "\pi_f"] \\
Z & X \ar[l, "f"]
\end{tikzcd}$$
is usually called a \emph{pullback diagram}.
Note that if $Z = \Jeden$ then $f = 1_X$, $g = 1_Y$, and the pullback is the product $X \times Y$.
We shall need the following well known basic properties of pullbacks.

\begin{lm}\label{LMcztyrivlstnsti}
	Let $\triple{W}{\pi_f}{\pi_g}$ be the pullback of $\pair fg$.
	\begin{enumerate}[itemsep=0pt]
		\item[{\rm(1)}] If $f$ is a surjection then so is $\pi_g$.
		\item[{\rm(2)}] If $X,Y,Z$ are topological spaces then $\pi_f, \pi_g$ are continuous, assuming $W$ has the topology inherited from the product $X \times Y$.
		\item[{\rm(3)}] For every nonempty set $S$, for every mappings $\map p S X$, $\map q S Y$ satisfying $f \cmp p = g \cmp q$ there exists a unique mapping $\map h S W$ such that $\pi_f \cmp h = p$, $\pi_g \cmp h = q$.
		\item[{\rm(4)}] If in {\rm(3)} the sets $Z,X,Y,S$ are topological spaces and $f,g, p, q$ are continuous, then $h$ is continuous.
	\end{enumerate}
\end{lm}

\begin{pf}
	(1) Given $y \in Y$, choose $x \in f^{-1}(g(y))$. Then $\pair xy \in W$ and $\pi_g(x,y) = y$.
	
	(2) Obvious.
	
	(3) If $h$ satisfies the assertion of (3) then necessarily $h(s) = \pair{p(s)}{q(s)}$ for $s \in S$.
	This is well defined, because $f(p(s)) = g(q(s))$.
	
	(4) The definition of $h$ clearly shows that it is continuous whenever $p,q$ are continuous.
\end{pf}

Clauses (2) and (4) of Lemma~\ref{LMcztyrivlstnsti} say, in the language of category theory, that the forgetful functor from the category of topological spaces to the category of sets creates pullbacks and, in particular, finite products. 

The following fact will be crucial for proving the main results. In the language of Boolean algebras, it can be found in Kechris \& Rosendal~\cite[Prop. 2.3]{KecRos}.

\begin{lm}\label{LMkrucjal}
	Assume $Z,X,Y$ are nonempty finite strictly positive probability spaces, $\map f X Z$, $\map g Y Z$ are measure preserving surjections.
	Let $\triple{W}{\pi_f}{\pi_g}$ be the pullback of $\pair fg$.
	Then there exists a strictly positive probability measure $\prob_{f,g}$ on $W$ such that $\pi_f$, $\pi_g$ become measure preserving.
	If the probabilities on $X,Y,Z$ are rational then so is $\prob_{f,g}$.
\end{lm}

\begin{pf}
	Fix $\pair xy \in W$ and let $z = f(x) = g(y)$.
	Define
	$$\prob_{f,g} \sn{\pair xy} = \frac{\prob \sn x \prob \sn y}{\prob \sn z}$$
	and extend $\prob_{f,g}$ to all subsets of $W$.
	Fix $x_0 \in X$ and let $z_0 = f(x_0)$. We have
	\begin{align*}
		\prob \pi_f^{-1}(x_0) &= \prob_{f,g} \setof{\pair {x_0}y}{g(y) = z_0} = \sum_{y \in Y,\; g(y)=z_0} \prob_{f,g} \sn{\pair{x_0}{y}} \\
		&= \sum_{y \in Y,\; g(y)=z_0} \frac{\prob \sn{x_0} \prob \sn y}{\prob \sn{z_0}} = \frac{\prob \sn{x_0}}{\prob \sn{z_0}} \cdot \sum_{y \in Y,\; g(y)=z_0} \prob \sn y \\
		&= \frac{\prob \sn{x_0}}{\prob \sn{z_0}} \prob g^{-1}(z_0) = \frac{\prob \sn{x_0}}{\prob \sn{z_0}} \prob \sn{z_0} = \prob \sn{x_0}.
	\end{align*}
	Hence, $\pi_f$ is measure preserving. In particular, $\prob_{f,g}$ is indeed a probability measure.
	By symmetry, $\pi_g$ is measure preserving too. Obviously, $\prob_{f,g}$ is rational, whenever the probabilities on $X,Y,Z$ are rational.
\end{pf}

Note that the lemma above remains valid when the set of values of the measures is replaced by any subset $V$ of the unit interval (possibly including the zero), as long as the following implication holds:
\begin{equation}
	\al \loe \beta < \gamma \in V \implies \frac{\al \beta}{\gamma} \in V.
\tag{$\divideontimes$}\label{EqAMmesrs}
\end{equation}
Note that the set $Q_2$ consisting of all dyadic rationals in the interval $[0,1]$ fails (\ref{EqAMmesrs}), as for example, we have
$$\frac13 = \frac{\frac 12 \cdot \frac 12}{\frac 34}.$$
On the other hand, if the measures of all points are of the form $1/2^{k}$, then the amalgamation is possible. We shall use this observation in the last section, extending the main results to different homogeneous probability measures.

\subsection{Boolean algebras with measures}

We now discuss some more details concerning the duality between profinite spaces with probability measures and Boolean algebras with finitely additive probability measures.
Namely, a \emph{probability measure} on a Boolean algebra $\Be$ is a function $\map \mu \Be {[0,1]}$ satisfying $\mu(1_\Be) = 1$ and $\mu(a \join b) = \mu(a) + \mu(b)$ for every $a, b \in \Be$ such that $a \meet b = 0_\Be$.
Here, $1_\Be$ is the unit of $\Be$ and $\join$ and $\meet$ denote the Boolean addition and multiplication operations, respectively.
The measure $\mu$ is \emph{strictly positive} if $\mu(a) = 0$ implies $a = 0_\Be$.
The Boolean algebra $\Be$ may be arbitrary, no assumptions on $\sig$-completeness are made, therefore this definition essentially differs from that of a countably additive probability measure on a complete Boolean algebra, see Fremlin's chapter in the Handbook of Boolean Algebras~\cite{HBAFremlin}. A measure in the sense of our definition is sometimes called a \emph{finitely additive} probability measure.
In any case, via Stone duality, Boolean algebras with probability measures correspond precisely to profinite spaces with Borel probability measures. Actually, we have defined profinite spaces as inverse limits of \emph{sequences} of finite sets, the more general (and more common) definition uses inverse systems, allowing arbitrary weight of the limit, and hence arbitrarily large cardinality of the Boolean algebra.

A measure preserving surjection between finite sets (or, more generally, between profinite spaces) corresponds to an embedding of the Boolean algebras preserving the associated measures. In particular, if $\Aa \subs \Be$ are two Boolean algebras with measures $\mu_\Aa$ and $\mu_\Be$ then this embedding is \emph{measure preserving} exactly if $\mu_\Be$ extends $\mu_\Aa$.
Now it is evident that every countable Boolean algebra with a probability measure is the union of a chain of finite subalgebras with the same (restricted) measure and the embeddings are obviously measure preserving. This, via Stone duality, provides straightforward proofs of Lemmas~\ref{LMdwaAjedn} and~\ref{LMciongimiarp}.
Furthermore, Lemma~\ref{LMomegaSmlns}, when transfered to Boolean algebras, becomes merely trivial.
Proposition~\ref{PROPsdbAIFHU} says that every countable Boolean algebra admits a strictly positive probability measure (this may easily fail for arbitrary Boolean algebras). Here, the proof transfers to showing that if $\Aa \subs \Be$ are finite Boolean algebras then every strictly positive probability measure on $\Aa$ extends to a strictly positive probability measure on $\Be$. This is easily proved by examining the atoms of both algebras, which is actually the same as moving to the corresponding Stone spaces.

Summarizing, it is clear that this note could have been written entirely in the language of Boolean algebras with probability measures. We have decided to use the language of profinite spaces for two reasons: Firstly, measure theory originated in the classical setting of sets or topological spaces, measures seem to be more natural when their domains are algebras of sets instead of abstract Boolean algebras.
Secondly, Lemma~\ref{LMkrucjal} seems to be much easier to manage, especially when tested with different value sets of the measures, see Section~\ref{SecDiscsFinsRems}. This lemma deals with pullbacks, which have a very concrete set-theoretic definition using the product. In order to translate it to Boolean algebras, one needs to use pushouts, i.e., suitable quotients of free sums which are more abstract and perhaps less intuitive objects\footnote{Some algebraists will obviously disagree here.}, see~\cite[Prop. 2.3]{KecRos}.

\section{The universal homogeneous rational probability}\label{SecPfThmOneJdn}

Note that the category $\kfP$ of finite rational strictly positive probability spaces is an inverse \fra\ category.
Indeed, Lemma~\ref{LMkrucjal} shows that it has the inverse amalgamation property. Furthermore, it has a terminal object, therefore it is inverse directed.
Finally, it has countably many isomorphic types and its hom-sets are finite. The following fact is rather trivial.

\begin{prop}
	The pair $\kfP \subs \kP$ satisfies {\rm(L0)}--{\rm(L3)}.
\end{prop}

Let $\pair{U_\infty}\Prob$ be the (unique up to measure preserving homeomorphisms) $\kfP$-projective object in $\kP$, existing by Theorem~\ref{ThmFraisseReversedOne}.
We show below that $U_\infty$ is a Cantor set and the measure $\Prob$ is universal homogeneous, proving Theorems~\ref{THMmejnOne} and~\ref{THMmejnTwo}. In the next section we construct universal homogeneous embeddings, thus proving Theorems~\ref{THMmejnThree} and~\ref{THMmejnFour}.

\begin{pf}[Proof of Theorem~\ref{THMmejnOne}]
	We first check that $\pair{U_\infty}{\Prob}$ satisfies (\los). Fix a clopen set $U \subs U_\infty$ with $r = \Prob U > 0$.
	If $r<1$ then let $S = \dn 01$ be the two-element probability space with $\prob \sn 0 = r$.
	If $U = U_\infty$ and $r=1$, then let $S = \Jeden$, where $\Jeden = \sn0$ is the trivial probability space.
	Let $\map q{U_\infty}S$ be such that $q^{-1}(0) = U$. Then $q$ is continuous and measure preserving, i.e., $q \in \kP$.
	Assume $r = r_0+r_1$ with $r_0, r_1>0$ rational and let $\map f T S$ be a surjection such that $f^{-1}(0) = \dn {a_0}{a_1}$ and $f^{-1}(1) = \sn b$ (in case $S \ne \Jeden$).
	Here, $T = \{a_0, a_1, b\}$ in case $S \ne \Jeden$ and $T = \dn{a_0}{a_1}$ otherwise.
	Define a probability $\prob$ on $T$ by setting $\prob \sn {a_i} = r_i$ for $i<2$.
	Then $f$ becomes measure preserving, i.e., $f \in \kfP$.
	As $U_\infty$ is $\kfP$-projective, there is a continuous measure preserving map $\map h{U_\infty} T$ such that $f \cmp h = q$.
	Let $U_i = h^{-1}(a_i)$, $i<2$.
	Then $U_0, U_1$ are disjoint clopen sets, $\Prob U_i = \prob \sn{a_i} = r_i$ for $i<2$, and $U_0 \cup U_1 = \inv{h}{\dn{a_0}{a_1}} = \inv h{f^{-1}(0)} = q^{-1}(0) = U$.
	
	Condition (\los) implies that the space has no isolated points, therefore $U_\infty$ is homeomorphic to the Cantor set.
	
	In order to show uniqueness, it suffices to check that (\los) implies $\kfP$-projectivity.
	So let $V \in \obj{\kP}$ satisfy (\los) and fix a $\kP$-arrow $\map q V S$ with $S$ finite. Fix a $\kfP$-arrow $\map f T S$.
	Fix $s \in S$ and let $f^{-1}(s) = \{a_0, \dots, a_{n-1}\}$. Let $r_i = \prob \sn{a_i}$.
	Then $r_0 + \dots + r_{n-1} = \prob \sn s = \Prob U$, where $U = q^{-1}(s)$.
	Condition (\los) together with an obvious induction show that there is a clopen partition $\{U_0, \dots, U_{n-1}\}$ of $U$ such that $r_i = \Prob U_i$ for $i<n$.
	Define $h$ on $U$ so that $h^{-1}(a_i) = U_i$ for $i<n$.
	We do the same for each $s \in S$, thus defining $h$ on the whole of $V$. Clearly, $h$ is continuous, measure preserving, and $f \cmp h = q$.
	
	It remains to check that the product of all finite strictly positive rational probability spaces satisfies (\los). Let $\sett{\pair{F_n}{\prob_n}}{\ntr}$ be a fixed enumeration of all these spaces and let $\prob_\infty$ be the product measure on $F_\infty := \prod_{\ntr} F_n$.
	Fix a clopen set $U \subs F_\infty$ with $r := \prob_\infty U$ and assume $r_0 + r_1 = r$, where $r_0,r_1>0$ are rational.
	Let $\al_0, \al_1$ be such that $r_i = \al_i r$, $i=0,1$. Then $\al_0 + \al_1 = 1$.
	Note that $U = U' \times \prod_{n>n_0}F_n$, where $U' \subs \prod_{n \loe n_0}F_n$.
	Fix $k > \max_{n \loe n_0} |F_n|$.
	There exists $m$ such that $F_n$ has at least $k$ elements and can be decomposed into $A_0$, $A_1$ so that $\prob_m A_i = \al_i$. By the choice of $k$ we know that $m > n_0$.
	Consider
	$$U_i = U' \times A_i \times \prod_{n > n_0, n \ne m} F_n, \qquad i=0,1.$$
	Then $\prob_\infty U_i = r \prob_m A_i = r \al_i = r_i$, which shows (\los).
\end{pf}

We are now ready to give a proof of Theorem~\ref{THMmejnTwo} in the special case where $K = \emptyset$.
Fix a nonempty clopen set $U \subs \cantor$ and define
$$\Prob_U(A) = \frac{\Prob(A)}{\Prob(U)}$$
for every Borel set $A \subs U$.
Then $\Prob_U$ is a rational probability measure on $U$ and $\pair{U}{\Prob_U}$ satisfies (\los).
Indeed, if $V \subs U$ is a nonempty clopen set with $\Prob_U V = r_0+r_1$, where $r_0,r_1$ are positive rationals, then $\Prob V = r_0 \Prob U + r_1 \Prob U$ therefore, by $(\los)$ we find a clopen partition $\dn{V_0}{V_1}$ of $V$ such that $\Prob {V_i} = r_i \Prob U$ for $i=0,1$; finally $\Prob_U{V_i} = r_i$ for $i=0,1$.
Hence, $\pair U{\Prob_U}$ is isomorphic in $\kP$ to $\pair{\cantor}{\Prob}$.
It follows that if $U$, $V$ are clopen subsets of $\cantor$ with $\Prob U = \Prob V$ then there exists a measure preserving homeomorphism from $U$ onto $V$.
The assertion of Theorem~\ref{THMmejnTwo} (with $K = \emptyset$) follows immediately from this property, namely, the required homeomorphism is obtained by gluing together fixed measure preserving homeomorphisms $\map {h_i}{U_i}{V_i}$, $i < n$, which exist by the remarks above.

We note the following easy consequence of property (\los), needed in the next section.

\begin{lm}\label{LMsfdkbsdjbvas}
	Let $W \subs \cantor$ be a nonempty clopen set, let $A_0, A_1 \subs W$ be disjoint measure zero sets, let $r_0, r_1 > 0$ be rational numbers such that $\Prob W = r_0 + r_1$.
	Then there exists a clopen partition $W = W_0 \cup W_1$ such that $A_i \subs W_i$ and $\Prob W_i = r_i$ for $i=0,1$.
\end{lm}

\begin{pf}
	By the regularity of $\Prob$, there are disjoint clopen neighborhoods $U_0$, $U_1$ of $A_0$, $A_1$, respectively, such that $U_0 \cup U_1 \subs W$ and $\Prob U_i = \delta_i < r_i$ for $i=0,1$.
	Applying (\los) to $W \setminus (U_0 \cup U_1)$, we obtain disjoint clopen sets $V_0$, $V_1$ such that $V_0 \cup V_1 = W \setminus (U_0 \cup U_1)$ and $\Prob V_i = r_i - \delta_i$ for $i=0,1$.
	Finally, $W_i := V_i \cup U_i$, $i=0,1$, are as required.
\end{pf}

\section{Universal homogeneous embeddings}

Towards proving Theorems~\ref{THMmejnThree} and~\ref{THMmejnFour}, we will work in a certain comma category related to $\kP$. Comma categories in \fra\ theory have already been successfully used by Pech \& Pech~\cite{PechPech, PechPech2}.

Throughout this section we fix a profinite space $K$. In other words, $K$ is a nonempty compact 0-dimensional second countable topological space.
We define the category $\fL_K$ as follows.

The objects of $\fL_K$ are continuous mappings $\map f K X$, where $X \in \obj{\kP}$.
Given two $\fL_K$-objects $\map {f_0}K{X_0}$, $\map {f_1}K{X_1}$, an $\fL_K$-arrow from $f_1$ to $f_0$ is a continuous measure preserving surjection $\map q {X_1}{X_0}$ satisfying $q \cmp f_1 = f_0$, as shown in the diagram below.
$$\begin{tikzcd}
& & & & X_1 \ar[dd, two heads, "q"] \\
K \ar[rrrru, "f_1"] \ar[rrrrd, "f_0"'] & & & & \\
& & & & X_0
\end{tikzcd}$$
The composition is the usual composition of mappings.
We define $\fK_K$ to be the full subcategory of $\fL_K$ whose objects are those $\map f K X$ with $X$ finite.

\begin{prop}
	$\fK_K \subs \fL_K$ satisfy {\rm(L0)}--{\rm(L3)} and $\fK_K$ is an inverse \fra\ category.
\end{prop}

\begin{pf}
	The first part is rather obvious.
	Concerning the second part, only the inverse amalgamation property requires an argument.
	On the other hand, this is a straightforward consequence of Lemma~\ref{LMkrucjal} combined with the crucial property of the pullback, stated in (3), (4) of Lemma~\ref{LMcztyrivlstnsti}.
\end{pf}

Note that $\fK_K$ has a terminal object $\map t K \Jeden$, therefore condition (G2) implies (G1). In other words, a $\fK_K$-projective object is $\fK_K$-generic in $\fL_K$.

\begin{tw}\label{THMfdngonsog}
	Let the Cantor set $\cantor$ be endowed with the universal homogeneous rational probability measure $\Prob$. Let $\map \eta K \cantor$ be a continuous mapping such that $\img \eta K$ is of measure zero in $\cantor$.
	Then $\eta$ is $\fK_K$-generic in $\fL_K$.
\end{tw}

\begin{pf}
	It suffices to show that $\eta$ satisfies condition (G2), which translates to the following:
	\begin{enumerate}
		\item[(E)] Given finite rational probability spaces $X, Y \in \obj{\kfP}$, given a measure preserving surjection $\map f Y X$, given a continuous mapping $\map b K Y$, given a continuous measure preserving surjection $\map g \cantor X$, there exists a continuous measure preserving surjection $\map h \cantor Y$ such that
		\begin{equation}
			f \cmp h = g \oraz h \cmp \eta = b.
			\tag{$\star$}\label{EqBigStrr}
		\end{equation}
	\end{enumerate}
	$$\begin{tikzcd}
	K \ar[dd, "b"'] \ar[rr, "\eta"] & & \cantor \ar[dd, two heads, "g"'] \ar[lldd, two heads, dashed, blue, "h"] \\
	&&&& \\
	Y \ar[rr, two heads, "f"']& & X
	\end{tikzcd}$$
	We say that a surjection $\map f Y X$ is \emph{prime} if there is $x_0 \in X$ such that $|f^{-1}(x_0)| = 2$ and $|f^{-1}(x)| = 1$ for every $x \in X \setminus \sn{x_0}$. Note that every measure preserving surjection between finite probability spaces is a composition of prime measure preserving surjections. Thus, it is enough to show (E) in case where $f$ is prime.
	So, assume $f^{-1}(x_0) = \dn{y_0}{y_1}$ with $y_0 \ne y_1$ and $f$ is one-to-one on $Y \setminus \dn{y_0}{y_1}$.
	
	Let $W = g^{-1}(x_0)$.
	There is a unique way to define $h$ on $\cantor \setminus W$ so that~(\ref{EqBigStrr}) holds.
	It remains to define $h$ on $W$.

	Let $A_i = \img \eta {b^{-1}(y_i)}$, $i=0,1$.
	Then $A_0$, $A_1$ are closed measure zero subsets of $W$ (it may happen that $A_i = \emptyset$).
	Applying Lemma~\ref{LMsfdkbsdjbvas}, we get a clopen partition $W = W_0 \cup W_1$ such that $\Prob W_i = \prob_Y \sn{y_i}$ and $A_i \subs W_i$ for $i=0,1$, where $\prob_Y$ is the probability measure of $Y$.
	Define $h \rest W$ so that $h^{-1}(y_i) = W_i$, $i=0,1$.
	It is clear that $h$ is measure preserving and satisfies $f \cmp h = g$.
	We check that $h \cmp \eta = b$.
	Fix $t \in K$ such that $b(t) \in \dn{y_0}{y_1}$ (otherwise there is nothing to check). Assume $b(t) = y_i$. Then $\eta(t) \in A_i \subs W_i$ and hence $h(\eta(t)) = y_i = b(t)$. Thus (\ref{EqBigStrr}) holds.
	This completes the proof.
\end{pf}

We are now ready to prove our main results.

\begin{pf}[Proof of Theorem~\ref{THMmejnTwo}]
	The existence of the partition $\{U_0, \dots, U_{n-1}\}$ follows from Theorem~\ref{THMmejnOne}.
	Let $\eta$ denote the inclusion $K \subs \cantor$.
	Let $S = \{0, \dots, n-1\}$ and define the probability measure $\prob$ on $S$ by setting $\prob \sn i = \Prob U_i = \Prob V_i$.
	Define $\map a K S$ in such a way that $U_i \cap K = a^{-1}(i) = V_i \cap K$ for $i<n$.
	Then $a$ is a $\fK_K$-object.
	Let $\map p \cantor S$, $\map q \cantor S$ be given by $p^{-1}(i) = U_i$, $q^{-1}(i) = V_i$, $i < n$. Then $p$, $q$ are measure preserving continuous surjections satisfying $p \cmp \eta = a = q \cmp \eta$. In other words, $p$ and $q$ are two $\fL_K$-arrows from $\eta$ to $a$.
	By the projective homogeneity of $\eta$ (Theorem~\ref{ThmFraisseReversedTwo}(2)) there exists an $\fL_K$-isomorphism $\map h \cantor \cantor$ satisfying $q = p \cmp h$.
	In other words, $h$ is a measure preserving homeomorphism of $\cantor$ such that $h \cmp \eta = \eta$ and $\img h {U_i} = V_i$ for $i<n$. The equation $h \cmp \eta = \eta$ simply means $h \rest K = \id K$.
\end{pf}

\begin{pf}[Proof of Theorem~\ref{THMmejnThree}]
	Let $\map \eta K \cantor$ be $\fK_K$-generic and fix a rational strictly positive probability measure $P$ on $K$.
	The identity mapping $\map{\id K}{K}{K}$ is an object of $\fL_K$, therefore by the projective universality of $\eta$ (Theorem~\ref{ThmFraisseReversedTwo}) there exists a $\fL_K$-arrow $\map f \cantor K$ from $\eta$ to $\id K$. In other words, $f$ is a continuous measure preserving surjection satisfying $f \cmp \eta = \id K$.
\end{pf}

\begin{pf}[Proof of Theorem~\ref{THMmejnFour}]
	Let $\eta$ denote the inclusion $K \subs \cantor$ and let $\xi = h \cmp \eta$. By Theorem~\ref{THMfdngonsog}, both $\eta$ and $\xi$ are $\fK_K$-generic, therefore by uniqueness (Theorem~\ref{ThmFraisseReversedOne}) there exists an $\fL_K$-isomorphism $\map H \cantor \cantor$ from $\eta$ to $\xi$. In other words,
	$H$ is a measure preserving homeomorphism satisfying $H \cmp \eta = \xi$, and this equation translates to $H \rest K = h$.
\end{pf}

\section{Discussion and final remarks}\label{SecDiscsFinsRems}

One can consider probability measures whose values on clopen sets are in a fixed countable subset of $[0,1]$. We only need to be sure that Lemma~\ref{LMkrucjal} on pullbacks is still valid, so that the comma category has inverse amalgamations.
For example, we can drop the assumption that the measures are strictly positive, thus adding zero to the value set.
Lemma~\ref{LMkrucjal} remains true, because (using the notation from its proof) if $\prob{\sn z} = 0$, then $\prob{\sn x} = 0$ and $\prob{\sn y} = 0$, where $z = f(x) = g(y)$, and we define $\prob_{f,g}\sn{\pair x y} = 0$. 

So, fix a countable set $V \subs [0,1]$ with $1 \in V$ and consider profinite spaces with probability measures whose values on clopen sets are in $V$.
Such a measure $\Prob_V$ will be called \emph{$V$-homogeneous} if it satisfies the following condition.
\begin{enumerate}
	\item[(\losv{V})] For every nonempty clopen set $U$, for every $r_0,r_1 \in V$ with $r_0 + r_1 = \Prob_V U$, there exists a partition $\dn{U_0}{U_1}$ into nonempty clopen sets such that $\Prob_V U_i = r_i$ for $i=0,1$.
\end{enumerate} 
In particular, let $\Probzero = \Prob_V$, where $V = \Qyu \cap [0,1]$. Thus $\Probzero$ differs from $\Prob$ by allowing the zero as a possible value at nonempty clopen sets.
Clearly, $\Probzero$ is a probability measure on the Cantor set $\cantor$. Let $C$ be the support of $\Probzero$. Then $C$ is a closed nowhere dense subset of $\cantor$, homeomorphic to the Cantor set, and $\Probzero$ restricted to $C$ is isomorphic to $\Prob$. Theorem~\ref{THMfdngonsog} now has to be modified by imposing that the image of $K$ is nowhere dense and disjoint from $C$.
Theorem~\ref{THMmejnFour} now says that every homeomorphism between closed nowhere dense subsets of the Cantor set that are disjoint from $C$ extends to a measure preserving auto-homeomorphism of $\cantor$.

We now discuss other possibilities. Namely, fix a countable set $V \subs [0,1]$ with $1 \in V$ and suppose that a $V$-homogeneous measure $\Prob_V$ exists. Denote by $\bfK_V$ the profinite space on which $\Prob_V$ lives.
Let $\kPv V$ and $\kfPv V$ denote the corresponding categories of profinite and finite probability spaces where now the values at clopen sets are in $V$.
Note that the $V$-homogeneity of $\Prob_V$ implies that $\pair{\bfK_V}{\Prob_V}$ is $\kfPv V$-projective (see the proof of Theorem~\ref{THMmejnOne}).
In particular, the measure $\Prob_V$ is unique, up to a measure preserving homeomorphism.

Let us now assume that the set $V$ satisfies the following conditions:
\begin{enumerate}[itemsep=0pt]
	\item[(H0)] $1 \in V$.
	\item[(H1)] $\al, \beta \in V \implies |\beta - \al| \in V$.
\end{enumerate}
Note that the existence of a $V$-homogeneous measure does not imply (H1). Actually, we may have corrected condition (\losv{V}) requiring additionally that for each $r \in V$ there is a clopen set of measure $r$, but this would give us only a weaker variant of (H1), namely,
\begin{enumerate}
	\item[] $\al \in V \implies 1 - \al \in V$.
\end{enumerate}
Take $V = \{\al, 1 - \al, 1\}$, where $0 < \al < \frac 12$ and $\al \ne \frac 13$.
Then there are no $r_0, r_1 \in V$ satisfying $r_0 + r_1 = 1 - \al$, therefore the appropriate two-element probability space is $V$-homogeneous.
This indicates that (H1) is a rather natural condition.
Note the following consequence of (H0) and (H1):
\begin{enumerate}
	\item[(H2)] $\al, \beta \in V, \; \al + \beta < 1 \implies \al+\beta \in V$.
\end{enumerate}
Indeed, assuming $\al \loe \beta$ we get $1 - \beta \in V$ and $\al < 1 - \beta$, therefore $(1 - \beta) - \al \in V$ and hence also $\al+\beta \in V$.

Coming back to the discussion of $V$-homogeneous measures, let us note that if $0 \in V$ then $\bfK_V$ is the Cantor set and the support of $\Prob_V$ is nowhere dense. Furthermore, $\Prob_V$ restricted to its support is isomorphic to $\Prob_{V'}$, where $V' = V \setminus \sn 0$.
So let us assume that $0 \notin V$.

Suppose first that $V$ is finite and let $r = \min V$.
Then $r = \frac 1m$ for some positive integer $m$.
Indeed, let $m$ be such that $(m-1)r < 1 \loe mr$.
Then $r \goe \frac 1m$ and $1-(m-1)r \in V$, therefore $1-(m-1)r \goe r$, which gives $r \loe \frac 1m$.
It follows that $V = \{ \frac 1m, \frac 2m \dots, 1\}$.
Let $\Prob_m$ denote the uniform probability measure on the $m$-element set $m = \{0, \dots, m-1\}$.
It is clear that $\Prob_m$ is $V$-homogeneous.

Suppose now that $V$ is infinite.
By (H1) and (H2), $V$ is dense in $[0,1]$ and hence $\bfK_V$ is the Cantor set.
Thus, an obvious question arises: are the main results (in particular, Theorems~\ref{THMmejnThree} and~\ref{THMmejnFour}) still valid?
The answer is evidently affirmative if $V = F \cap (0,1]$, where $F$ is a countable subfield of $\Err$, as in this case condition (\ref{EqAMmesrs}) holds and all the other arguments can be repeated.
Furthermore, following the last part of the proof of Theorem~\ref{THMmejnOne}, we can easily deduce that $\Prob_m^\nat$, the countable infinite power of $\Prob_m$, is $V$-homogeneous, where $V$ is the set of all positive $m$-adic rationals in the unit interval.
Recall that $r \in \Err$ is an \emph{$m$-adic rational} if it has a finite expansion with base $m$ (in case $m=2$, such numbers are called \emph{dyadic rationals}).
Now, Lemma~\ref{LMkrucjal} is not true (see the example after its proof), however it becomes true when the measure on $Z$ has atoms of the form $\frac{1}{m^i}$. Thus we can restrict our category $\kfPv V$ by requiring that the measures have values of the form $\frac{1}{m^i}$ at each point (atom). By this way, we still obtain a version of Theorem~\ref{THMmejnFour}, however it is not clear how to obtain Theorem~\ref{THMmejnThree}.

Of course, we could also consider the extreme case of reducing the value set $V$ to $\sn0$, therefore forgetting the measures at all. By this way we would obtain the original Knaster-Reichbach result saying that every homeomorphism between closed nowhere dense subsets of the Cantor set extends to an auto-homeomorphism.

In general, we do not know how to characterize countable sets $V \subs [0,1]$ for which there exists a $V$-homogeneous probability measure.

\separator

Finally, different relational/algebraic structures on profinite spaces should lead to similar results, however most likely it will be difficult to characterize the embeddings that are \fra\ limits in the corresponding comma category.
A typical example here could be the category of finite groups. It has pullbacks, therefore it clearly has the inverse amalgamation property. The inverse \fra\ limit is the universal projectively homogeneous second countable profinite group. It is not clear, however, how the corresponding homogeneity results concerning closed nowhere dense subgroups should look like.

\paragraph{Acknowledgments.}
The authors would like to thank the anonymous referee for extremely careful reading of the manuscript, pointing out several improvements, in particular, for simplifying the diagram in the proof of Theorem~\ref{THMfdngonsog}.
The authors are also grateful to Grzegorz Plebanek for providing reference~\cite{BaSu} and to Karol Baron for useful comments concerning measure theory.
The second author would like to thank the Erwin Schr\"odinger International Institute for Mathematics and Physics (Vienna, Austria) where part of the research was being conducted (December 2016).
The research of all the authors was supported by grant No. 16-34860L awarded jointly by the Austrian Science Foundation (FWF) and the Czech Science Foundation (GA\v CR).

\end{document}